\documentclass[sn-mathphys-num]{sn-jnl}

\usepackage{graphicx}%
\usepackage{multirow}%
\usepackage{amsmath,amssymb,amsfonts}%
\usepackage{amsthm}%
\usepackage{mathrsfs}%
\usepackage[title]{appendix}%
\usepackage{xcolor}%
\usepackage{textcomp}%
\usepackage{manyfoot}%
\usepackage{booktabs}%
\usepackage{algorithm}%
\usepackage{algorithmicx}%
\usepackage{algpseudocode}%
\usepackage{listings}%
\usepackage{tikz}
\usepackage{tkz-euclide}

\theoremstyle{thmstyleone}%
\theoremstyle{thmstyletwo}%

\theoremstyle{thmstylethree}%
\newtheorem{definition}{Definition}%

\begin{document}

\title[Generating maps on oriented surfaces]{Generating maps on oriented surfaces using the homomorphism principle}
\author*[1]{\fnm{Gunnar} \sur{ Brinkmann}} \email{Gunnar.Brinkmann@UGent.be}
\affil*[1]{\orgdiv{Toegepaste Wiskunde, Informatica en Statistiek}, \orgname{Universiteit Gent}, \orgaddress{\street{Krijgslaan 281 S9}, \city{Gent}, \postcode{B 9000}, \country{Belgium}}}

\abstract{
  
  In this article we describe an algorithm that can be applied for the generation
  of various classes of maps on orientable surfaces. It uses existing generators for abstract graphs and combines them with an efficient embedding and isomorphism rejection routine.
  The generation rate of the program implementing the algorithm depends a lot on the class of
  maps to be generated, but is quite high -- more than a million non-isomorphic structures per second -- for some relevant classes of maps.
  The same program can also be used to embed specific graphs on a given orientable surface in all non-isomorphic ways. It  can  serve as a tool in many applications
  where classes of maps on orientable surfaces are studied and provides a very general independent test for specialized generation programs. We also give enumeration results
  for 3-regular, 4-regular, and 5-regular maps as well as all maps and some maps with just one face.
  }

\keywords{graph, genus, structure enumeration, homomorphism principle}
\pacs[MSC Classification]{05C10, 05-04, 05C85}

\maketitle

\section{Introduction}

Polyhedra -- or from a combinatorial point of view: 3-connected maps of genus 0 -- belong to the first structures where computers were used for enumeration \cite{G65}.
Later many algorithms and programs for the enumeration of various classes of polyhedra and (fewer) algorithms and programs for  the generation of
maps of higher genus were developed and used in research to test
conjectures and find suitable structures. See e.g. \cite{genusclique}\cite{genusdependence}. The program plantri \cite{plantri} implements algorithms for several classes of plane maps -- among others
triangulations, quadrangulations, and bipartite maps. Though plantri looks like one program, the generation of the classes is based on very different construction methods, in most cases
given by a recursive structure of the class, so that plantri is
more a conglomerate of several programs than a single program.
In \cite{Sulanke} Sulanke describes the algorithm that is the basis for his program {\em surftri} and gives results obtained with it
for seven orientable and non-orientable surfaces. Surftri can generate all triangulations of these surfaces and by manipulation of
the triangulations also e.g. all maps with a given number of vertices on those surfaces. Nevertheless it would not be possible to efficiently generate maps with certain properties like e.g.\
4-regular maps or bipartite maps. This would -- just like in plantri -- need specialized operations and isomorphism rejection routines.
The number of papers where structure enumeration is applied to study plane maps exceeds the
number of papers  where structure enumeration is applied for the study of maps of higher genus by far. This is also due to the fact that for most classes of maps on nontrivial genus
there are simply no generation programs available. To this end more generation programs for maps on surfaces of higher genus should be developed.

In \cite{4regplane} C-L.~Simon and B.~Stucky use plantri to generate 4-regular graphs from which they constructed {\em pinning semi-lattices}. They asked
whether a similar tool for the generation of 4-regular maps on surfaces of higher genus -- especially the torus -- does exist. Of course this is a very natural
question as next to 3-regular maps 4-regular maps are probably the most interesting class, as such maps can also serve as models for knots. In this article we describe an
easy algorithm based on \cite{genuscomp} that can e.g.\ generate 3-regular, 4-regular, and 5-regular maps, bipartite maps, or maps with a given degree sequence.

Not distinguishing between mirror images, the total number of embeddings of a graph $G$ with maximum degree at least $3$ and
vertices $\{1,\dots ,n\}$ is $\frac{1}{2}\prod_{i=1}^n (\deg(i)-1)!$ with $\deg(i)$ the degree of vertex $i$.
For a single $4$-regular graph with trivial symmetry this means $3\cdot 6^{n-1}$ non-isomorphic embeddings in at most $\frac{n+3}{2}$ different genera -- while it is known
that if it is 3-connected, at most one is in genus $0$. Already this simple observation 
makes clear that even for $4$-regular graphs -- let alone graphs with larger degrees -- due to the enormous growth it will only be possible and useful to generate complete  list of maps
for relatively small vertex numbers and genera. So the methods must focus on these achievable numbers for which generation is possible and not on asymptotic behaviour of the
algorithms. Using the described method to generate $r$-regular maps of genus $g$ with $n$ vertices is inefficient if the fraction of $r$-regular graphs that can be embedded in an orientable surface of genus $g$
is small. Asymptotically for each $r\ge 3$ and $k$, $r$-regular graphs almost surely contain a $K_k$ minor \cite{minor_reg}, so for each fixed $r$ and $g$ the ratio of regular graphs that can be embedded in genus
$g$ goes to zero, but as the tables will show, in most cases this is long after the number of structures to be generated (and in applications: to be tested) becomes too large for any possible use of the structures.

It seems to be naive and inefficient to generate classes of maps on surfaces by generating possible underlying abstract graphs and embedding them on the
  surface in every possible way -- especially as it is often not clear which graphs can be embedded in the given surface.  While this is surely the case for the plane,
  we will show that with carefully
  implemented efficient algorithms for embedding the graphs and isomorphism rejection of the maps, this method can be sufficiently fast for some cases
  and in others possibly even close to the
  best possible approach. The efficiency of the algorithm described depends strongly on the class to be generated. It varies from several million
  non-isomorphic maps per second in cases where each input graph has many embeddings to almost none for classes where only a tiny fraction of the
  input graphs allows an embedding in the surface. Fortunately many interesting classes fall into the first case.
  The algorithm will also be applicable for the generation of maps with a prescribed number of faces instead of a prescribed genus -- e.g.\ all maps with $7$ vertices and
  only one face, embedded in an orientable surface.

\section{Basic Definitions}

All graphs considered in this article are simple, undirected, and connected.  A map -- also called a combinatorially embedded graph -- is a triple $M=(V,E,s)$
with $G=(V,E)$ the {\em underlying abstract graph} in
which the undirected edges $\{u,v\}\in E$ are considered to consist of two oppositely oriented edges $[u,v]$ incident to $u$ and $[u,v]^{-1}=[v,u]$ incident to $v$.
For each vertex of a map, the oriented edges incident to it are assigned a cyclic ordering which can be considered the clockwise order around the common vertex when the graph is drawn.
We denote the set of oriented edges as $E_o$ and for a given oriented edge $e$ the function $s():E_o\to E_o$ maps $e$ onto the next oriented edge in the cyclic order.
This is also called a {\em rotation system}.
An isomorphism between maps $(V,E,s)$ and $(V',E',s')$ is a graph isomorphism $\phi:V\to V'$ of the underlying graphs that respects the cyclic order at each vertex (in
which case it is called orientation preserving) 
or reverses the cyclic order at each vertex (in which case it is called orientation reversing).
Or in exact words (with the induced action of $\phi$ on $E_o$): we have that $\phi(s(e))=s'(\phi(e))$ for all $e\in E_o$ or
$\phi(s(e))=(s')^{-1}(\phi(e))$  for all $e\in E_o$.
Isomorphism classes of these maps have a
one-to-one correspondence to the topologically defined homeomorphism classes of cellular embeddings of simple graphs
\cite{top_graph_theory}\cite{graphs_on_surfaces}. An angle in a map is a pair $e,e'$ of oriented edges with $e'=s(e^{-1})$.
A face is a cyclic sequence of oriented edges, so that every two edges $e,e'$ following each other in this sequence form an angle. The genus of a map is defined as
$g=1-(|V|-|E|+|F|)/2$ with $F$ the set of faces of the map and it is the same as the genus of the corresponding topological embedding.

For isomorphism rejection, we will apply the {\em homomorphism principle}. The homomorphism principle
is often used in structure enumeration -- see e.g. \cite{LBFL80},\cite{GKL92} or \cite{BCH01}, but the rest
of the article can also be understood without prior knowledge about the underlying principles of the homomorphism principle.

\section{The Algorithm}.

The basic algorithm to generate a certain class of maps on a given oriented surface works as follows:

\begin{description}
\item[a.)] Generate a set of pairwise non-isomorphic abstract graphs that contain all underlying graphs of the maps in the class to be generated.
  For this task in most cases suitable algorithms and programs already exist.

\item[b.)] Assign cyclic orders around each vertex to embed the abstract graphs in all ways in the prescriped genus avoiding or removing isomorphic copies.
\end{description}

Already at this point it is visible how the homomorphism principle can be exploited: each isomorphism between two maps induces an isomorphism between the
underlying graphs -- or in this case we can even say: {\bf is} an isomorphism between the underlying graphs. As in this algorithm isomorphic maps can only be
generated by assigning different orientations to the same graph, this isomorphism must be an automorphism. For underlying graphs with a trivial automorphism
group this means that no isomorphism rejection is necessary: the only possible isomorphism is the identity, so different cyclic orderings of the oriented edges
also imply non-isomorphic maps and no isomorphism rejection is necessary.

For graphs with nontrivial automorphisms, only the -- often very few -- automorphisms of the graph must be tested as isomorphisms of the maps.

\begin{figure}[tb]
	\centering
        \resizebox{0.65\textwidth}{0.65\textwidth}
        {
\begin{tikzpicture}[scale=0.06]
\def\vertexscale{1.20}
\def\labelscale{1.60}
\node [circle,black,draw,scale=\vertexscale] (1) at (0.48113,-1.22089) {1};
\node [circle,black,draw,scale=\vertexscale] (2) at (48.44187,-54.20977) {2};
\node [circle,black,draw,scale=\vertexscale] (3) at (-47.07447,40.43962) {3};
\node [circle,black,draw,scale=\vertexscale] (4) at (33.42866,-30.73762) {4};
\node [circle,black,draw,scale=\vertexscale] (5) at (-21.23975,-46.52692) {5};
\node [circle,black,draw,scale=\vertexscale] (6) at (44.89099,25.60256) {6};
\node [circle,black,draw,scale=\vertexscale] (7) at (-49.53796,-8.20160) {7};
\node [circle,black,draw,scale=\vertexscale] (8) at (-2.02685,46.37409) {8};
\tkzDefPoint(-89.10065,45.39905){9}
\tkzDefPoint(-70.71068,-70.71068){10}
\tkzDefPoint(-98.76883,15.64345){11}
\tkzDefPoint(-98.76883,-15.64345){12}
\tkzDefPoint(-89.10065,-45.39905){13}
\tkzDefPoint(-50.00000,-86.60254){14}
\tkzDefPoint(-25.88190,-96.59258){15}
\tkzDefPoint(0.00000,-100.00000){16}
\tkzDefPoint(25.88190,-96.59258){17}
\tkzDefPoint(50.00000,-86.60254){18}
\tkzDefPoint(89.10065,45.39905){19}
\tkzDefPoint(98.76883,15.64345){20}
\tkzDefPoint(98.76883,-15.64345){21}
\tkzDefPoint(89.10065,-45.39905){22}
\tkzDefPoint(-50.00000,86.60254){23}
\tkzDefPoint(-25.88190,96.59258){24}
\tkzDefPoint(-0.00000,100.00000){25}
\tkzDefPoint(25.88190,96.59258){26}
\tkzDefPoint(50.00000,86.60254){27}
\tkzDefPoint(-70.71068,70.71068){28}
\tkzDefPoint(70.71068,70.71068){29}
\tkzDefPoint(70.71068,-70.71068){30}
\draw [black] (1) to (8);
\draw [black] (1) to (6);
\draw [black] (1) to (4);
\draw [black] (1) to (5);
\draw [black] (1) to (7);
\draw [black] (2) to (4);
\draw [black] (2) to (22);
\node [draw=none,fill=none,scale=\labelscale] () at (93.55569,-47.66900) {7};
\draw [black] (2) to (18);
\node [draw=none,fill=none,scale=\labelscale] () at (52.50000,-90.93267) {6};
\draw [black] (2) to (17);
\node [draw=none,fill=none,scale=\labelscale] () at (27.17600,-101.42221) {8};
\draw [black] (3) to (23);
\node [draw=none,fill=none,scale=\labelscale] () at (-52.50000,90.93267) {5};
\draw [black] (3) to (8);
\draw [black] (3) to (7);
\draw [black] (3) to (9);
\node [draw=none,fill=none,scale=\labelscale] () at (-93.55569,47.66900) {6};
\draw [black] (4) to (16);
\node [draw=none,fill=none,scale=\labelscale] () at (0.00000,-105.00000) {8};
\draw [black] (4) to (5);
\draw [black] (4) to (6);
\draw [black] (4) to (21);
\node [draw=none,fill=none,scale=\labelscale] () at (103.70728,-16.42562) {7};
\draw [black] (5) to (15);
\node [draw=none,fill=none,scale=\labelscale] () at (-27.17600,-101.42221) {8};
\draw [black] (5) to (14);
\node [draw=none,fill=none,scale=\labelscale] () at (-52.50000,-90.93267) {3};
\draw [black] (5) to (7);
\draw [black] (6) to (27);
\node [draw=none,fill=none,scale=\labelscale] () at (52.50000,90.93267) {2};
\draw [black] (6) to (19);
\node [draw=none,fill=none,scale=\labelscale] () at (93.55569,47.66900) {3};
\draw [black] (6) to (20);
\node [draw=none,fill=none,scale=\labelscale] () at (103.70728,16.42562) {7};
\draw [black] (6) to (8);
\draw [black] (7) to (13);
\node [draw=none,fill=none,scale=\labelscale] () at (-93.55569,-47.66900) {2};
\draw [black] (7) to (12);
\node [draw=none,fill=none,scale=\labelscale] () at (-103.70728,-16.42562) {4};
\draw [black] (7) to (11);
\node [draw=none,fill=none,scale=\labelscale] () at (-103.70728,16.42562) {6};
\draw [black] (7) to (8);
\draw [black] (8) to (25);
\node [draw=none,fill=none,scale=\labelscale] () at (-0.00000,105.00000) {4};
\draw [black] (8) to (26);
\node [draw=none,fill=none,scale=\labelscale] () at (27.17600,101.42221) {2};
\draw [black] (8) to (24);
\node [draw=none,fill=none,scale=\labelscale] () at (-27.17600,101.42221) {5};
\tkzDefPoint(-68.55786,72.79986){A}
\tkzDefPoint(68.55786,72.79986){B}
\tkzDefPoint(0.0,0.0){C}
\tkzDrawArc[<-,line width=0.9mm, gray](C,B)(A)
\tkzDefPoint(72.79986,68.55786){A}
\tkzDefPoint(72.79986,-68.55786){B}
\tkzDefPoint(0.0,0.0){C}
\tkzDrawArc[<-,line width=0.9mm, gray](C,B)(A)
\tkzDefPoint(68.55786,-72.79986){A}
\tkzDefPoint(-68.55786,-72.79986){B}
\tkzDefPoint(0.0,0.0){C}
\tkzDrawArc[->,line width=0.9mm, gray](C,B)(A)
\tkzDefPoint(-72.79986,-68.55786){A}
\tkzDefPoint(-72.79986,68.55786){B}
\tkzDefPoint(0.0,0.0){C}
\tkzDrawArc[->,line width=0.9mm, gray](C,B)(A)
\node [black,circle,draw,fill=white,scale=0.75,line width=1mm] (10) at (-70.71068,-70.71068) {};
\node [black,circle,draw,fill=white,scale=0.75,line width=1mm] (28) at (-70.71068,70.71068) {};
\node [black,circle,draw,fill=white,scale=0.75,line width=1mm] (29) at (70.71068,70.71068) {};
\node [black,circle,draw,fill=white,scale=0.75,line width=1mm] (30) at (70.71068,-70.71068) {};
\end{tikzpicture}

          }
	\caption{A graph with $8$ vertices and nontrivial symmetry that has the largest number of embeddings in the orientable surface of genus $3$ of all such graphs.
          The graph is drawn with a minimum genus embedding on the torus where as usual opposite sides have to be identified.
          The only nontrivial automorphism exchanges the two vertices with degree $7$ (that is: vertex $7$ and vertex $8$) and stabilizes the others.
          In an orientable surface of genus $1$ it has $49$ non-isomorphic embeddings, in genus $2$ it has
        $611\ 031$ non-isomorphic embeddings, in genus $3$ it has $494\ 136\ 137$ non-isomorphic embeddings, and in genus $4$ it has $70\ 375\ 033\ 157$ non-isomorphic embeddings.}
	\label{fig:manyembeddings}
\end{figure}
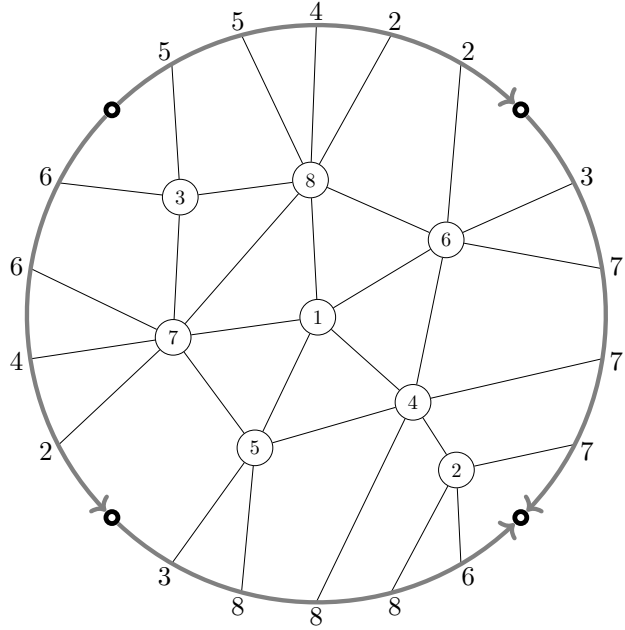

One could avoid the output of isomorphic maps by maintaining a list that is newly initialized for every input graph and contains {\em canonical representatives} for
each map generated. When using lists, the definition of the canonical representative does not even have to have any connection with the construction.
This isomorphism rejection with lists is in fact also part of the implementation, but only
for testing purposes. Unfortunately already some small graphs with nontrivial automorphisms -- so graphs where the lists would have to be used --
have extremely many non-isomorphic embeddings already for relatively small genera, so that the lists would become very large and inefficient and could even cause
memory problems.
An example is the graph with nontrivial automorphism group in Figure~\ref{fig:manyembeddings} that has only $8$ vertices, but already $70\ 375\ 033\ 157$ non-isomorphic embeddings in genus $4$.

\medskip

We work with a canonical
representative of each class:

\begin{definition}\label{def:aut}
  Let $M=(V,E,s)$ with $V=\{1,\dots ,n\}$ be a map. For a vertex $v$ with smallest neighbour $w_1$ the string $s_M(v)$ is
  defined as $w_1 w_2 \dots  w_{\deg(v)}$ with for $2\le i \le \deg(v)$ the vertex $w_i$ given by $[v,w_i]=s([v,w_{i-1}])$.
  The string $S_p(M)$ (with ``p'' for {\em orientation preserving}) is defined as $s_M(1)  s_M(2)  \dots  s_M(n)$.
  The string $S_r(M)$ (with ``r'' for {\em orientation reversing}) is defined as $s^{-1}_M(1)  s^{-1}_M(2)  \dots  s^{-1}_M(n)$ with  $s^{-1}(v)$ defined analogously to $s_M(v)$,
  but with $[v,w_i]=s^{-1}([v,w_{i-1}])$ for $2\le i \le \deg(v)$.

  The string $S(M)$ is the lexicographically smallest one of $S_p(M)$ and $S_r(M)$.

  The automorphism group of a graph $G$ or map $M$ is denoted by $Aut(G)$, resp.\ $Aut(M)$.
  An automorphism $\phi: V \to V$ of the underlying graph $G=(V,E)$ induces an action on the oriented edges and we define the image $\phi(M)=(V,E,s')$
  by $s'(e)= \phi(s(\phi^{-1}(e)))$.

  For a map  $M=(V,E,s)$ we define the canonical string for the isomorphism class (with fixed underlying graph $(V,E)$) as
  $S_{can}(M)=\min \{ S(\phi(M)) \mid \phi \in Aut((V,E))\}$.

\end{definition}

\medskip

For a map $M=(V,E,s)$ and an automorphism $\phi$ of $(V,E)$ the map $\phi(M)=(V',E',s')$ is isomorphic to $M$ by the automorphism $\phi$:

As $\phi$ is an automorphism of $(V,E)$, it is an isomorphism between the identical underlying graphs $(V,E)$ and $(V',E')$.
Furthermore for $e\in E$ we have  \\ $s'(\phi(e))= \phi(s(\phi^{-1}(\phi(e))))=\phi(s(e))$.

If on the other hand $M=(V,E,s)$ is isomorphic to $M'=(V,E,s')$ by a w.l.o.g.\ orientation preserving isomorphism $\psi$ (otherwise choose $(s')^{-1}$),
for each oriented edge $e'$ in $M$ we have
$s'(\psi(e'))=\psi(s(e'))$. Applying this for an arbitrary oriented edge $e\in E_o$ to $e'=\psi^{-1}(e)$ we get $s'(\psi(\psi^{-1}(e)))=\psi(s(\psi^{-1}(e)))$, so
$s'(e)=\psi(s(\psi^{-1}(e)))$ and therefore $s'$ is exactly of the form described in Definition~\ref{def:aut} and for $M=(V,E,s)$ the set $\{\phi(M) | \phi \in Aut((V,E)\}$
contains only isomorphic maps and at least one  member for each set of isomorphic maps with the same underlying graph.

\medskip

\begin{definition}\label{def:canon}

  A map $M=(V,E,s)$ is called the canonical representative of its class (with respect to the underlying graph $G=(V,E)$), if $S(M)=S_{can}(M)$.

\end{definition}

\medskip

So maps where the underlying graph has a trivial automorphism group are always canonical representatives of their class. Maps with a nontrivial automorphism group
have in general also embeddings that are not canonical representatives of their class as some automorphisms of the graph are not automorphisms of the map.
An exception are of course  3-connected plane maps -- due to Whitney's theorem.

The algorithm described in \cite{genuscomp} can compute all maps with a given genus for a given input graph $G=(V,E)$. This algorithm fixes the orientation around one vertex,
so it constructs exactly one of
$(V,E,s)$ and $(V,E,s^{-1})$ for each rotation system that gives a map of the required genus.

In our algorithm a map is accepted if and only if it is the canonical representative of its class, so that no lists have to be maintained and memory consumption due to large lists of maps is not a problem.
This is essentially the principle of orderly generation as described by Read \cite{Read78}.
As soon as there is a vertex of degree at least $3$,
there are two canonical representatives of each class, as mirror images have the same representative  --  but only one of them is generated by the construction algorithm.

The test for being the canonical representative can be applied to each map $(V,E,s)$ that is constructed: for each automorphism $\phi$ of $(V,E)$ the string
$S(\phi(M))$ is constructed and compared to $S(M)$. Then $M$ would be accepted if and only if no smaller string is found. In fact for some cases this turned out to be
the fastest method. We will discuss the reason later on, but will also describe a method that is in general much faster.

For most cases it was -- as expected -- faster to try to detect non-canonicity already early: the construction routine described in  \cite{genuscomp}
starts with a subgraph that can (up to inversion of the cyclic order) be uniquely embedded and then inserts the remaining edges one by one. The order in which the edges are embedded
is according to a predefined plan. For an edge $e$ we call the iteration in which it is planned to be inserted $t(e)$.
In exceptional cases some edges can be inserted earlier than planned, but at iteration $t(e)$, at least all edges $e'$ with $t(e')<t(e)$ are already embedded. To this end we know
for each iteration, for which automorphisms $\phi$ of $(V,E)$ a prefix of $S(\phi(M))$ -- with $M$ a possible map constructed from this partial embedding -- can already be computed.
For $n'\le n$ the parts $s_{\phi(M)}(1)  s_{\phi(M)}(2)  \dots  s_{\phi(M)}(n')$ and $s^{-1}_{\phi(M)}(1)  s^{-1}_{\phi(M)}(2)  \dots  s^{-1}_{\phi(M)}(n')$ can be computed as soon as
all edges incident to $\phi(1), \dots ,\phi(n')$ are embedded. If also all edges incident to $1, \dots ,n'$ are embedded, the prefixes of $S(M)$ and $S(\phi(M))$ can already be compared.
This can have a positive impact in several ways:

\begin{itemize}
  
\item If the prefix of $S(\phi(M))$ is smaller, we can stop the construction and backtrack.

\item  If the prefix of $S(\phi(M))$ is larger, the automorphism $\phi$ will not be considered for any test of a map constructed from this one by adding further edges.

\item Even if the prefix of $S(\phi(M))$ is identical to that of $S(M)$, there is an advantage: the prefix tested is only compared once -- and not for each map constructed
  from the partial embedding. Once the prefix up to vertex $n'$ is compared and found identical, for all successors, the test starts at the neighbours of $n'+1$.

\end{itemize}

A similar distribution of the partial canonicity tests over the recursion tree was already a key principle in \cite{minibaum}. Although this method can have a large positive impact,
in this case there are also problems:

\begin{itemize}
  
\item One needs the automorphism group already before the construction starts and not only after the first map is constructed. In cases where the group is not part of the input,
  it has to be computed even for input graphs that do not have a single embedding in the required genus -- which is not the case if the group is only applied for complete embeddings.

\item In order to apply the method early, edges incident to vertex $1$ and to images of vertex $1$ must be inserted early. This is in conflict with the order
  chosen in \cite{genuscomp} that focuses on small branching and chooses edges for which -- a priori -- only a small number of possible ways to insert them
  exist. E.g. for optimal efficiency of the construction, in a $4$-regular graph, edges between two vertices that have degree $2$ at the moment are preferred to
  edges with two endpoints of degree $3$. For optimal efficiency of the bounding via isomorphisms it would be the other way around -- especially if the
  vertices of degree $3$ are the image of vertices with a small label under an automorphism.

\end{itemize}

To this end, testing for isomorphisms only at the end can be faster -- especially if the ratio of input graphs that allow an embedding with the chosen genus is
small and graphs that allow such an embedding have a relatively small number of them in average. In such cases finding embeddings or deciding that there are
none is the most difficult and time consuming part of the algorithm. An example is the generation of $4$-regular maps with $15$ vertices.  In order to focus on
the essential differences, the following relative times are only for graphs with a non-trivial automorphism group and the time for generating the input graphs is not
counted. In that case, for genus $1$ testing canonicity only when all edges are embedded is fastest, but even testing isomorphisms with lists is only slightly slower
and more than two times faster than intermediate testing. This difference is almost entirely due to the different order in which edges are inserted --
computing the automorphism group for fewer graphs is only a very small fraction of the difference.  For genus $2$, embedding is easier and isomorphism
rejection plays a bigger role: for genus $2$, intermediate testing is $1.15$ times faster than testing only at the end and $1.42$ times faster than using lists.
In these cases also a compromise -- using the order that is optimal for the construction, but nevertheless doing intermediate canonicity checks -- was slower than
testing only at the end.

The program uses some heuristics to decide whether to do intermediate tests. In cases where the user has a good reason to assume that one of the methods would surely perform better, it
can be chosen as an option when starting the program. If the genus of the surface on which the graphs must be embedded is e.g.\ close to the genus of the graphs,
not doing intermediate tests is often better.

In order to have an insertion order that performs well for the construction as well as for isomorphism rejection, non-regular graphs are sometimes relabeled in a way
that vertices with small degree get small labels. Nevertheless, vertex $1$ is -- if possible -- always a vertex with degree at least $3$, as this way it can be determined
whether $s$ or $s^{-1}$ give the smaller string as soon as all edges of vertex $1$ are inserted.

The algorithm relies on other generation programs for generating the underlying abstract graphs. The default is to compute the automorphism group by using {\em nauty}
\cite{nauty2}. By default nauty computes generators for the automorphism group, but the package also contains functions computing each element of the group.
The program described here uses the full group. In general this is
a very bad idea as the groups can become extremely large. Nevertheless, for the applications that this
algorithm is intended for, the groups are at least in average relatively small. In some cases, the room needed for storing the group can be a problem though and although very general, the program can
not be applied for all applications -- it would e.g.\ not be a suitable algorithm to prove that the star with 20 vertices has a unique embedding in the plane.
Often the generation programs have some idea about the automorphism groups and e.g. the program {\em
  minibaum} \cite{minibaum} generating cubic graphs in fact knows -- and can output -- the whole automorphism group together with the graph constructed.
In cases where the computation of the group
takes an important part of the generation time, this gives a considerable speedup. For cubic maps with $20$ vertices and genus $1$ using intermediate testing
(which is the fastest option for cubic graphs) the speedup is more than $20\%$.

\section{Results}

We will now give some tables with numbers of maps and generation rates for different genera and graph classes. The first numbers in the second line of each block are the structures
generated (but not outputted) per second on an Intel i7-9700 restricted to 3 GHz running Linux with several jobs running at the same time. The numbers are
multiples of $1.000$ maps. The second numbers are the percentages of input graphs allowing at least one embedding.
The times for generating the graphs and computing the automorphism groups by minibaum (cubic), {\em genreg} \cite{MM98} and nauty (quartic, 5-regular and 6-regular), resp.\ {\em geng} \cite{McK96}
and nauty are included when computing the generation rates.
In the third, resp.\ third and fourth, line, the average number of embeddings per graph (rounded) and the maximum number of embeddings of a graph are given.
In cases where the running times are too small to give a meaningful generation rate, we write {\em nm}. As can be seen in the tables, the generation rates differ
a lot. They are very small when the required genus is relatively small for the class of graphs, so that
many graphs have very few or even no embeddings in that genus. E.g.:\ for $4$-regular maps of genus $1$ with $19$ vertices, only $0.52\%$ of the input
graphs allow an embedding in genus $1$ and only $395$ maps per secon are generated. Some additions to the program generating the input graphs that avoid at least many of the graphs with high genus would
help, but as the intention here is to have a very general tool, this was not implemented. On the other hand, for some classes of maps the generation rate is very high
-- even for the same class of underlying graphs, like $4$-regular maps of higher genus. For genus $3$ or $4$, the generation rate is sometimes more than $3$ or even $4$ million non-isomorphic $4$-regular
maps per second. Here we
have the opposite effect: up to $13$ vertices all input graphs even allow an embedding with genus $2$ and as a consequence many embeddings in the surfaces of genus $3$ and $4$.

Not only graphs with a small genus (or equivalently: many faces) are interesting. Also maps with few faces are relevant and well studied.
An interesting example is \cite{tomo14} where the self assembly of
polypeptide nanostructures is approached via {\em double traces}, which correspond to maps (orientable or non-orientable) with just one face. Note that
for orientable maps the Euler characteristic is always even, so that for a given number of vertices and edges the number of faces is either always even
or always odd, so that e.g.\ cubic orientable maps with $8$ vertices and just one face do not exist. Numbers of maps with just one face are given in Table~\ref{tab:1face}.

\section{Testing}

\begin{table}[h]
\begin{tabular}{c|c||c|c||c|c}
  $|V|$ & cubic & $|V|$ & quartic & $|V|$ & all\\
  \hline
  4        & 4,f: (i),(ii),(iii),(iv)  &  5  &  4,f: (i),(ii),(iii),(iv)   & 4   &   4,f: (i),(ii),(iii),(iv)   \\
  \hline
 6 &    4,f:   (i),(ii),(iii),(iv)  &  6  &  4,f: (i),(ii),(iii),(iv)  & 5  &   4,f: (i),(ii),(iii),(iv)   \\
  \hline
 8 &   4,f: (i),(ii),(iii),(iv)   &   7 &  4,f: (i),(ii),(iii),(iv) & 6  &   4,f: (i),(ii),(iii),(iv)   \\
  \hline
  10 &   4,f:  (i),(ii),(iii),(iv)  &   8 &  4,f: (i),(ii),(iii),(iv)  & 7  &   3: (i),(ii),(iii)   \\
   &    &    &    &   &   4: (i),(ii)   \\
  \hline
  12 &   4,f:  (i),(ii),(iii),(iv)   &   9 & 4,f: (i),(ii),(iii),(iv)  & 8  &  3: (i),(ii)  \\
  \hline
  14 &    4,f:  (i),(ii),(iii),(iv)   &   10 & 3: (i),(ii),(iii),(iv)  & 9  &  2: (i),(ii) \\
  \hline
  16&   4,f:  (i),(ii),(iii)  &    11&    2: (i),(ii),(iii),(iv)   & 10  &  1: (i),(ii)  \\
  \hline
 18  & 4,f:  (i),(ii),(iii)     &   12 &      4: (i),(ii)            &  &    \\
  \hline
 20 &  4: (i),(ii),(iii)     &    13 &     4: (i),(ii)             & &    \\
  \hline
 22  &  3: (i),(ii)  &    &    &    &    \\
\end{tabular}
\caption{ The cases for which the numbers of maps were compared. An example for the interpretation:  3,f: (i),(ii),(iii) stands for: the numbers of maps generated by methods  (i),(ii),(iii) were compared for genus 0 to 3 and for prescribed face numbers.
  The face numbers were always 1 and 2. }\label{tab:test}
\end{table}

Even if an algorithm is correct, it is still necessary to carefully test the implementation. Small implementation errors like typos can occur in parts of the code that are not often used, so that they
are not easily detected and independent tests should be performed as far as possible.

As mentioned before, not many enumerative results exist for surfaces of nontrivial genus. Although the program is not well suited for genus 0, comparing
existing results for genus $0$ is also a test for the embedding and isomorphism rejection routines. We compared the numbers of all plane maps up to 11 vertices,
all connected cubic plane maps up to $26$ vertices, and all 3-connected quartic plane maps on up to $17$ vertices to the
results of plantri and had complete agreement. Also surftri can generate classes of maps that can be independently generated by the program described here and
surftri also works for nontrivial genus. Unfortunately the numbers of structures generated by surftri differ already relatively early from ours: for all maps on the torus,
there is agreement up to 7 vertices, but the approach
described here has more than $2 000$ maps more than given in \cite{Sulanke}. We tested the maps for isomorphism and contacted the
author of surftri. It turned out that some graphs are missing in the lists that correspond to the numbers of all toroidal maps reported in \cite{Sulanke}.
In the meantime the error in the implementation has been detected and Thom Sulanke released a new version of surftri. Corrected numbers are given in Table~\ref{tab:all}.

In \cite{phdThomas} the results of a specialized program generating
cubic connected maps of the torus are presented. There was agreement up to $22$ vertices, but for $24$ vertices our number was $6$ more than reported in \cite{phdThomas}.
Also in this case it turned out that the 6 graphs were missing in the list of \cite{phdThomas}.

So for higher genus, some new and as far as possible independent tests had to be done.

The construction routine from \cite{genuscomp} has been seriously tested as described in  \cite{genuscomp}, so the main effort was on testing the isomorphism rejection routines.
Nevertheless due to small modifications and the interaction of the construction routines with the
isomorphism rejection routines, at least some independent tests also for that part were performed.

The program itself has three ways to reject isomorphic copies:

\begin{description}
\item[(i):] isomorphism rejection by canonicity tests of completely embedded graphs
\item[(ii):] isomorphism rejection by canonicity tests already during the construction
  \item[(iii):] isomorphism rejection by lists
\end{description}

The methods (i) and (iii) do not interact with the construction routines described in \cite{genuscomp}.

The most independent test was to compare the numbers of maps generated by the program to the ones obtained by generating
{\bf all} embeddings by constructing all cyclic orders around a vertex, filtering them for the genus (resp.\  the number of faces) that is wanted  and then removing isomorphic copies by a routine
using lists and the canonicity routines from plantri. This is of course an extremely slow method and due to the large number of maps only possible for small cases. Let us call
this exhaustive method (iv).

The number of non-isomorphic structures produced by the various methods were compared as described in Table \ref{tab:test} and there was complete agreement..

\section{Conclusions}

Structure enumeration programs are applied to test conjectures, find examples, or get an intuition on the behaviour of some invariants. In this article an algorithm is presented that can be applied in
a very broad context, as for many classes of abstract -- that is: non-embedded -- graphs, generation programs already exist. The generation rate varies a lot depending on the class, but for several
useful classes, the generation will not be the bottleneck when non-trivial computations are performed on the generated structures. When only a tiny part of the possible input graphs allows an
embedding in the chosen surface, this approach is not suitable -- e.g.\ for generating triangulations of surfaces and using all graphs with the suitable number of edges (e.g. $|E|=3|V|$ for the torus)
as input.

As soon as the genus is two or even three or four, the number of maps also grows very fast for most classes, so that in many cases even programs developed especially for one of the possible classes
would allow at most one or two vertices more -- even if no expensive tests are applied to the output.
Developing specialized algorithms might nevertheless still be interesting for some classes of graphs on the torus or when it is not possible to generate sets of input graphs where a large
fraction can be embedded. Also for graph classes with very large automorphism groups special programs would still be useful.
For these algorithms, the present program may still turn out to be a useful test of the implementations of the algorithms.

As a final conclusion one can say that the described program can turn out to be very useful for many different applications when maps on orientable surfaces are studied.
The program described in this article as well as the programs generating the input graphs can be obtained from the author, resp.\ from the web page \cite{nauty2}.

\begin{table}[h]
\begin{tabular}{r|r|r|r|r}
  $|V|$ & genus 1 & genus 2 & genus 3 & genus 4 \\
  \hline
  4 & 2 &  0  & 0   &  0  \\
  &   nm ; 100\%      &    &    &    \\
  &   2 ; 2      &    &    &    \\
  \hline
  6 & 7 & 3 &  0  &  0  \\
  &   nm ; 100\%      &     nm ; 100\%      &    &    \\
  &   3.5 ; 5  & 1.5 ; 2   &    &    \\
  \hline
  8 & 37    &  71 & 0   & 0   \\
  &  nm ; 100\%    &   nm ; 100\%     &    &    \\
  &   7.4 ; 15 & 14.2 ; 28   &    &    \\
  \hline
  10 &  232  & 1 234 &  437    &  0  \\
  &   nm ; 100\%      &     nm ; 100\%      &   nm ; 94.7\%  &    \\
  &   12.2 ; 27  &  65 ; 180  & 23 ; 72   &    \\
  \hline
  12 &  1 742 &  20 087   & 30 096    &  0  \\
  &    nm ; 100\%    &   nm ; 100\%     &  nm ; 100\%     &    \\
  &   20.5 ; 79  & 236.3 ; 864 & 354 ;  1 328   &    \\
  \hline
  14 &  14 379   &  305 671    &  1 185 527  &  347 879   \\
  &    nm ; 98.04\%   &    nm ; 100\%     & 1 852   ; 100\%    &  633 ; 96.3\%  \\
  &   28.3 ; 150 &  600.5 ; 2 024  &  2 329 ; 5 480  &  683 ; 1 176  \\
  \hline
  16 &  127 282  &  4 394 543    &  34 513 779  &  39 835 633   \\
  &   nm ; 94.83\%   &   1.400; 100\%    & 2 271 ; 100\%   &  1 772 ; 99.9\% \\
  &   31.4  ; 260   & 1 082 ; 4 408   & 8 501 ;  17 616  & 9 812 ; 21 440   \\
  \hline
  18 &  1 170 626 & 60 297 999  &  828 743 496  &   2 372 429 287  \\
  &    275 ; 87.14\%    &    1 166; 100\%    & 2 130 ; 100\%  &  2 326 ;  100\% \\
  &   28.3 ; 560  &  1 460 ; 11 272  & 20 065 ;  55 392   &  57 442  ; 86 880\\
  \hline
  20 &   11 021 519 &   797 468 502 & 17 414 556 400   &   99 107 654 171   \\
  &    181 ; 73.83\%   &    874; 100\%    & 1 833 ; 100\%  & 2 387   ;  100\%  \\
  &    21.6 ; 1 248 &   1 562 ;  28 432 & 34 113 ; 164 608   &  194 143 ;  310 272 \\
  \hline
  22 &  105 245 468   &   10 242 995 020  &  331 891 605 521  &    \\
  &    105 ; 56.46\%   &   610; 99.96\%    & 1 442 ; 100\%   &    \\
  &   14.4 ;  2 624  &  1 399 ; 64 800 & 45 344 ; 458 368   &    \\
  \hline
  24 &   1 014 547 272   & 128 521 948 069   &    &    \\
  &    55 ; 38.86\%   & 402 ; 99.7\%  &    &    \\
  &   8.6 ; 5 632      & 1 090  ; 151 360 &    &    \\
  \hline
  26 & 9 846 233 168    &    &    &    \\
       & 26 ;  24.17\%        &    &    &    \\
       &   4.7 ;   11 520     &    &    &    \\
 
\end{tabular}
\caption{ The number of cubic maps for different genera. Input are all connected cubic graphs with the given number of vertices.  }
\label{tab:cubic}
\end{table}

\begin{table}[h]
\begin{tabular}{r|r|r|r|r}
  $|V|$ & genus 1 & genus 2 & genus 3 & genus 4 \\
  \hline
  \hline
  6 & 2 & 1 &  0  & 0  \\
  &  nm ; 100\%   &  nm ; 100\% &  &   \\
  & 2 ; 2 & 1 ; 1 &  &   \\
   \hline
  8 & 5 & 8 &  0 & 0  \\
  &   nm ; 100\%  &   nm ; 100\%  &  &   \\
  &  5 ; 5 & 8 ; 8 &  &   \\
   \hline
  10 & 8  & 39 &  20 &  0 \\
  &   nm ; 100\% & nm ; 100\% & nm ; 100\% &   \\
  & 4 ; 6 & 19.5 ; 24 &  10 ; 14 &   \\
   \hline
  12 &  24  & 266  & 472 &  0 \\
  &   nm ; 100\% & nm ; 100\% &  nm ; 100\% &   \\
  & 4.8 ; 9 & 53.2 ; 102 & 94.4 ; 196  &  \\
   \hline
  14 & 66 & 1 605  &  7 732 &  2 624 \\
  & nm ;  84.6\% &  nm ; 100\%  &  nm ; 100\% &   nm ; 100\%  \\
  &  5.08 ; 12 & 123.5 ; 507 & 595 ; 1 390 & 202 ; 510   \\
   \hline
  16 & 208 &  9 447 & 99 279  &  138 692  \\
  & nm ; 79\% &   nm ; 100\%  & nm ; 100\% &  nm ; 100\%  \\
  & 5.47 ; 24  & 248.6 ; 816 & 2 613 ; 7 080 &  3 650 ; 10 036 \\
   \hline
  18 & 644 & 53 006  & 1 054 193  & 3 829 528  \\
  &  nm ; 63\% &   nm ; 100\%  & 994 ; 100\% & 1 272  ; 100\%\\
  &  4.46 ; 54 & 356 ; 1 654 & 7 075 ; 25 936 &  25 701 ;  85 632  \\
   \hline
  20 & 2 389 & 290 656 &  9 708 739 &  73 287 765 \\
  & nm ; 46.5\%  & 363 ; 100\% &  1 108 ; 100\% &  1 647 ; 100\% \\
  & 3.4 ; 56 & 413 ; 2820 & 13 810 ; 50 000 & 104 250 ; 248 224   \\
   \hline
  22 &  8 603  & 1 542 772  &  80 225 358  &  1 096 894 958 \\
  & 12.7 ; 29.1\% & 266 ; 99\% &  975 ; 100\% &  1 622 ; 100\% \\
  & 2.08 ; 110 & 373 ; 4 104 &  19 416 ; 92 144 &  265 463  ; 693 184 \\
   \hline
  24 &  32 816  & 8 012 818  & 609 867 387  &  13 743 014 491  \\
  & 5.7 ; 16.5\% & 157 ; 97.4\% &  744 ; 100\%& 1 433 ; 100\%  \\
  & 1.1 ; 214 & 271 ; 7 368 & 20 618 ; 173 424 &  464 621 ; 1 628 224 \\
   \hline
  26 &  126 059   & 40 757 091 &  4 339 316 832 &   \\
  & 2.1 ; 8.1\% & 81 ; 92.2\%  & 511 ; 100\% &   \\
  & 0.51 ; 436 & 166 ; 15 360 & 17 666 ; 337 408 &   \\
    \hline
  28 &  495 462  & 2 291 589 &  &   \\
  &  0.69 ; 3.7\% & 0.407 ; 79.4\% &  &   \\
  &  0.22 ; 888 & 89 ; 32 000 &  &   \\
     \hline
  30 &   1 958 864  &  &  &   \\
  & 0.2 ; 1.5\% &  &  &   \\
  & 0.08 ; 1 808 &  &  &   \\
  
\end{tabular}
\caption{ The number of cubic bipartite maps for different genera. Input are all connected cubic bipartite graphs with the given number of vertices.}
\label{tab:cubicbip}
\end{table}

\begin{table}[h]
\begin{tabular}{r|r|r|r|r}
  $|V|$ & genus 1 & genus 2 & genus 3 & genus 4 \\
  \hline
  5  & 6   &   31 & 13   &  0  \\
  &   nm ; 100\%  &   nm ; 100\%  &   nm ; 100\%    &    \\
  &  6 ; 6       & 31 ; 31 &  13 ; 13  &    \\
  \hline
  6 &  17 &  206 &  371   &  0  \\
  &   nm ; 100\%   &  nm ; 100\%    & nm ; 100\%  &    \\
  &  17 ; 17 & 206 ; 206 & 371 ; 371    &    \\
  \hline
  7 & 38  &  1 415   &  8 778    &   3 231   \\
  &     nm ; 100\%  &  nm ; 100\%    & nm ; 100\%  &   nm ; 100\%   \\
  &   19 ; 28  & 708 ; 1 066 &   4 389 ; 6 772  &   1 616 ; 2 457  \\
  \hline
  8 &  130   & 10 386     & 150 539 &   246 373   \\
  &    nm ; 100\%  &   nm ; 100\%    & 138 ; 100\%   &  123 ; 100\%  \\
  &    21.7 ; 50 & 1 731 ; 5 138  & 25 090 ; 77 392 &  41 062 ;  128 776 \\
  \hline
 9 &   441 &  69 727  & 2 041 804 &   9 217 119  \\
  &    nm ; 100\%  &   nm ; 100\% & 2 041 ; 100\%  & 2 275 ; 100\%  \\
 &  27.6 ; 91  & 4 358 ; 13 720 &  127 612 ;   374 487  &  576 070 ;  1 635 181\\
 \hline
 10 &   1 575  &  438 220   &   23 079 153  &  228 628 213    \\
  &  nm ; 98.3\% & 1 184  ; 100\%  & 2 487 ; 100\%  &  3056 ; 100\%   \\
 &  26.7 ; 120 &  7 427 ; 29 253 & 391 172 ;  &    3 875 054 ; \\
 &  &   &  1 376 303   &      12 744 766  \\
 \hline
 11 & 5 776  &  2 588 850    &  226 525 348  &   4 284 624 611   \\
  &  nm ; 91\% & 1 488 ; 100\%   & 3 146 ; 100\%   & 4 072 ; 100\%  \\
 &  21.8 ; 200 & 9 769 ; 39 082 & 854 813 ;  &  1 6168 395 ; \\
  &  &  2 404 712  & &37 118 207\\
 \hline
 12 &  22 147  &   14 585 067  & 1 990 911 841    &  65 358 458 068   \\
  &  nm ; 80\%  &  1 458 ; 100\% &  3 280 ; 100\% &  4 431 ; 100\%  \\
 &  14.3 ; 331 & 9 446 ; 81 987   & 1 289 451 ;   &  42 330 608 ; \\
  &   &   &  5 526 400   & 107 237 838 \\
 \hline
 13 & 85 211  &  78 993 084   &  16 021 555 075   &  850 836 552 595  \\
  & 85 ; 59.3\% & 1 161 ; 100\%  & 2 956 ; 100\%  & 4 266  ; 100\% \\
 &  7.9 ; 605 & 7 329 ; 120 977 & 1 486 505 ;    & 78 941 970 ;   238 140 746 \\
 &  &  & 9 085 101  &    \\
 \hline
 14 & 334 073  & 414 664 509   & 120 017 599 160   &    \\
  &  47 ; 36.9\% & 856 ; 100\%  & 2 589 ; 100\%  &    \\
 &  3.8 ;  678 & 4 703 ;  175 724  &   1 361 237 ;   &    \\
 &  &   &   14 874 995 &    \\
 \hline
 15 &  1 313 202  & 2 121 221 883    &    &    \\
  &  22 ; 19.3\% &  594 ;  99.97\% &    &    \\
 &  1.63 ; 1 312 & 2 633 ;  330 254 &    &    \\
 \hline
 16 & 5 198 341 & 10 624 052 933   &    &    \\
  &  9 ; 9\%       & 372 ;  99.6\%  &    &    \\
 &  0.65 ; 2 437 & 1 321 ;  644 440   &    &    \\
 \hline
 17 & 20 617 249 & 52 275 605 556   &    &    \\
  &  3.4 ; 3.8\%  & 215 ; 96.9\%   &    &    \\
 &  0.24 ; 4 994 & 606 ;   1 283 550  &    &    \\
 \hline
 18 &  82 024 953 &    &    &    \\
  &  1.2 ; 1.46\% &    &    &    \\
 &  0.08 ;  9 856 &    &    &    \\
 \hline
 19 & 326 807 983 &    &    &    \\
  & 0.395 ; 0.52\% &    &    &    \\
 & 0.027 ; 16 572  &    &    &    \\
 \end{tabular}
\caption{ The number of quartic maps for different genera.  Input are all connected quartic graphs with the given number of vertices.}\label{tab:quartic}
\end{table}

\begin{table}[h]
\begin{tabular}{r|r|r|r|r}
  $|V|$ & genus 1 & genus 2 & genus 3 & genus 4 \\
  \hline
  6 & 4  &  492  &   17 482   &   87 274   \\
  &   nm ; 100\%    &  nm ; 100\%  &   nm ; 100\%   &  2.4 ; 100\%  \\
  &  4 ; 4  & 492 ; 492  &   17 482 ;   17 482  &  87 274 ; 87 274   \\
  \hline

  8 &  18    &  21 577  &  5 072 241  &   227 148 023   \\
  &   nm ; 100\%     &   nm ; 100\%    &  453 ; 100\%   &  596 ; 100\%  \\
  &  6 ; 10   &  7 192 ; 15 807 &  1 690 747 &  75 716 008  \\
  &     &   &   3 682 089 &  163 897 234  \\
  \hline

 10 &  121    &  850 795 &  838 558 593  &  171 462 130 936  \\
  &  nm ; 51.7\%   &  567 ; 98.3\% & 1 480 ; 100\% & 1 372  ; 100\% \\
 &   2.02 ; 9 &  14 180 ;  51 607 &  13 975 977  &  2 857 702 182  \\
  &    &   & 50 094 526    & 9 892 860 202   \\
 \hline

 12 &  894    &  22 617 536  &   67 962 981 676  &    \\
  & 1.4 ; 4.9\%   & 203 ; 99.6\%  & 1 746 ; 100\%   &    \\
 & 0.11 ; 26    & 2 882 ; 41 867   &  8 659 911 ; 60 213 171 &    \\
 \hline

 14 &  7 141   & 464 461 076   &    &    \\
  &  0.035 ; 0.1\%   &  21.4 ; 86.5\% &    &    \\
 & $<$0.01 ; 55    & 134 ;  97 280 &    &    \\
\hline
16 &   61 628   &   &    &    \\
  & 0.00052 ;  0.0013\%   &  &    &    \\
 & $<$0.01 ; 130  &  &    &    \\

 \end{tabular}
\caption{ The number of 5-regular maps for different genera.  Input are all connected 5-regular graphs with the given number of vertices.}\label{tab:5-reg}
\end{table}

\begin{table}[h]
\begin{tabular}{r|r|r|r|r}
  $|V|$ & genus 1 & genus 2 & genus 3 & genus 4 \\
  \hline
 7 &  1 & 363   &  372 638    &   59 029 357   \\
  & nm ; 100\%    &  0.182 ; 100\%  & 0.218 ; 100\%   &  0.244 ; 100\%  \\
  & 1 ; 1  & 363 ; 363  & 372 638 ;    &  59 029 357 ;  \\
 &    &   &  372 638   &  59 029 357   \\
 
 8 & 1  & 1 285   &   3 477 656   &  1 475 568 462    \\
  &  nm ; 100\%    & 11.2 ; 100\%  &  28.6 ;100\%   & 37.7 ; 100\%  \\
  & 1 ; 1  & 1 285 ;  &  3 477 656 ;  &  1 475 568 462  ;  \\
 &    &  1285  &   3 477 656  &   1 475 568 462   \\
 
 9 & 2  &  5 804  &  43 507 193    & 44 977 095 936     \\
  & nm ; 50\%    &  46 ; 100\%  &  328 ; 100\%   & 504 ; 100\%  \\
  & 0.5 ; 1   & 1 451 ;  & 10 876 798 ;    &  11 244 273 984 ;  \\
 &    & 3 942   &  29 862 753   &  30 618 992 354  \\
 
 10 &  1 &  21 186  &   395 408 132   &  921 922 123 291    \\
  & nm ;  4.76\%    & 30 ; 100\%  &  1 038 ; 100\%   &  2 823 ; 100\%  \\
  &  1 ; 1  &  1 009 ; &  18 828 959 ;   &  43 901 053 490 ; \\
   &    & 4515  & 74 249 627    &  161 746 523 382  \\
 \end{tabular}
\caption{ The number of 6-regular maps for different genera.  Input are all connected 6-regular graphs with the given number of vertices.}\label{tab:6-reg}
\end{table}

\begin{table}[h]
\begin{tabular}{r|r|r|r|r|r}
  $|V|$ & genus 0 & genus 1 & genus 2 & genus 3 & genus 4 \\
  \hline
  4  & 6  &   3  & 0  &  0 &   0 \\
   & nm ; 100\%  &  nm ; 33.3\%     &    &   &    \\
   & 1 ; 1  & 0.5 ; 2      &   &   &    \\
   \hline
  5 & 25  & 70   & 81     &  13   &  0  \\
   & nm ; 95.2\%  &  nm ; 62\%  & nm ; 19.1\%   & nm ; 4.8\% &    \\
   &  1.2 ; 2 & 3.3 ; 17  & 3.9 ; 39  & 0.62 ; 13  &    \\
   \hline
  6 &  179   &  2 656  & 26 167  &  125 370   &    181 067  \\
   & nm ; 88.4\%  &  nm ; 82.1\%     & nm ; 46.4\%   & 14.8 ; 16.1\% &  4.5 ; 3.6\%  \\
   &  1.6 ; 6 &  23.7 ; 185 &  234 ; 4 182  &   1 119  ;  45 143  &  1 617 ;  87 274\\
   \hline
  7 & 2 014  &  126 466   & 5 423 933   &  150 911 857   &  2 341 220 550  \\
   & nm ; 75.7\%  &  nm ; 94.4\%  & 1 267 ; 74.3\%  & 72 ; 42.8\% & 8 ; 16.9\%   \\
  & 2.4 ; 24  &  148.3 ; 1 790  & 6 359 ; 159 231   &  11 814 074 &  2 744 690 \\
  &   &    &   &   11 814 074 &   347 476 067  \\
   \hline
  8 &  31 178 &  6 072 971   &  779 566 207   &   73 813 568 429 &    \\
   & nm ; 53.7\%  &  1 886 ; 98.7\%    & 4 283  ; 92.4\%  & 2 979 ; 74.6\%&    \\
  & 2.8 ; 80   & 546.2 ; 10 740  & 70 124 & 6 639 702 &    \\
   &    &   & 2 527 472 & 669 269 163 &    \\
   \hline
 9 & 580 555   &  280 313 481 &  86 650 642 553   &   &    \\
   & 231 ; 27.5\%  & 2 233 ; 98.2\%     & 4 578 ; 98.7\%  &  &    \\ 
 & 2.2 ; 240  &  1 074 ; 63 480 & 331 893  &   &    \\
 &  &   &   19 060 920 &   &    \\
  \hline
 10 & 12 046 072   & 12 392 228 369 &   &   &    \\
  & 121 ; 9\%  & 2 172 ; 87.7\%     &   &   &    \\
 & 1.03 ; 1080  &  1 058 ;  317 400    &   &   &    \\
 \hline
11 & 267 836 680 &  527 780 778 675 &   &   &    \\
  & 39 ; 1.73\%  &   1 841 ; 55.4\%   &   &   &    \\
  & 0.27 ; 3 780   & 524 ;  1 522 920    &   &   &    \\

 \end{tabular}
\caption{ The number of all maps for different genera.  Input are all connected graphs with the given number of vertices. Though the program is not meant and in general not useful for genus 0,
the numbers for genus 0 are mentioned in this case as an example, as in exceptional cases with low connectivity it might be useful.}\label{tab:all}
\end{table}

\begin{table}[h]
\begin{tabular}{r|r||r|r||r|r}
  $|V|$ & 3-regular & $|V|$ & 4-regular & $|V|$ & all \\
    \hline
    6  &        3                          & 5  &         13                          & 3  & 1 \\
    &       nm ;        100\%        &   &     nm ;          100\%         &   & nm ; 50\% \\
          &     1.5 ; 2                    &   &       13 ; 13                        &   & 0.5 ; 1\\
      \hline
      10 &    437                      &  7 &     3 231                           & 4  &  3 \\
      &    nm ;       94.74\%      & &       nm ; 100\%                   &   & nm ; 50\%\\
            &       23 ; 72               &  &   1615.5  ; 2457               &   & 0.5 ; 1 \\
        \hline
        14 &  14 379         &  9 &     2 925 883          & 5  & 33  \\
        &   nm ;   98.04\% &   &     770 ; 100\%          &   & nm ; 52.4\% \\
              &      28.25 ; 150    &   &    182 868 ;    535 725  &   & 1.57 ; 13 \\
          \hline
          18 &  581 195 456  &  11 &  4 113 902 484  & 6  &  47 953\\
          &    670 ;       97.48\%       &  & 1 178 ; 100\%         &   & 1.8 ; 50\%\\
                &     14 072 ;   31 744  & & 15 524 160 ; 34 939 998   &   & 428 ; 28 928\\
            \hline
              &                                 &    &                                       & 7  &  6 016 954 644 \\
                  &                                  &   &                                  &   & 316 ; 49.5\% \\
            &                                  &   &                                  &   & 7 053 874\\
                           &                                  &   &                                  &   & 5 133 871 944\\
              \hline

 \end{tabular}
\caption{ Results for some classes of maps with just one face. The graph with the maximum number of non-isomorphic embeddings with one face and 7 vertices is $K_7$ minus an edge. It has minimum genus 1 and embeddings with one face are in genus 7. }\label{tab:1face}
\end{table}



\end{document}